  \documentclass[authoryear,preprint,review,12pt]{elsarticle}

\usepackage{lineno}
\usepackage{amsfonts}
\usepackage{amssymb}
\usepackage{epsfig}
\usepackage{color}
\usepackage{amsmath}
\usepackage[normalem]{ulem}
\usepackage{epstopdf}
\usepackage{subfig}

\def\R{{\mathbb{R}}}
\def\C{\mathbb{C}}

\def\f{\varphi}

\journal{Chemical Engineering Science}

\begin{document}

\begin{frontmatter}

\title{Robust optimization of periodically operated nonlinear uncertain processes}

\author[label1]{Darya Kastsian}
\ead{darya.kastsian@rub.de}
\author[label1]{Martin M\"{o}nnigmann\corref{cor1}}
\ead{martin.moennigmann@rub.de}
\address[label1]{Automatic Control and Systems Theory, Ruhr-Universit\"{a}t Bochum, 44801 Bochum, Germany}
\cortext[cor1]{Corresponding author. Tel.: +49 234 3224060; fax: +49 234 3214155.}

\begin{abstract} 
We present a method for determining optimal modes of operation for autonomously oscillating systems with uncertain parameters. In a typical application of the method, a nonlinear dynamical system is optimized with respect to an economic objective function with nonlinear programming methods, and stability is guaranteed for all points in a robustness region around the optimal point. The stability constraints are implemented by imposing a lower bound on the distance between the optimal point and all stability boundaries in its vicinity, where stability boundaries are described with notions from bifurcation theory. We derive the required constraints for a general class of periodically operated processes and show how these bounds can be integrated into standard nonlinear programming methods. We present results of the optimization of two chemical reaction systems for illustration. 
\end{abstract}

\end{frontmatter}

\linenumbers

\section{Introduction}\label{sec:Intro} 
The impact of autonomous oscillations and periodic forcing on economic process performance has been investigated for decades.
For example, \cite{Douglas1966} demonstrate that the performance of an isothermal continuous stirred-tank reactor (CSTR) may be improved by periodic forcing of the feed. The authors also consider a first order irreversible exothermic reaction in a nonisothermal CSTR. For this case, they show that autonomous oscillations may lead to increased average product concentration compared to steady state operation.
Similar investigations have been carried out later by other authors. 
\cite{Jianquiang2000} use autonomous oscillations to improve the performance of a bioreactor used for 
sludge water treatment.
\cite{Stowers2009} show that oscillations can increase the product yield in yeast fermentation.
\cite{Parulekar2003} demonstrate that the performance of series-parallel reactions can be improved by forced periodic operation. The authors also discuss the benefit of forced periodic operation compared to steady state operation in recombinant cell culture processes.
\cite{Abashar2010} show that periodically forced fermentors provide higher average bioethanol concentrations than fermentors operated in steady state.

Whenever models of the production process of interest and its economics are available, it is an option to use linear or nonlinear programming methods to find an optimal mode of operation.
It is known, however, that optimizing a dynamical system in this way may result in a steady state or periodic mode of operation that, while optimal with respect to the economic objective, is unstable \citep{Moennigmann2002}. In general, optimal but unstable solutions are not useful in practice.

Approaches inspired by applied bifurcation theory have been used to state constraints on stability properties in optimization problems. Since these methods are based on normal vectors to manifolds of critical points such as bifurcations points, they are jointly referred to as the normal vector approach for short. 
Originally, the normal vector approach was developed to guarantee stability of optimal equilibria of ordinary differential equations (ODE) and differential-algebraic (DAE) systems \citep{Moennigmann2002,Moennigmann2007}. It has been applied to a number of examples from chemical engineering \citep{Moennigmann2003,Moennigmann2005}.   
\cite{Gerhard2008} and \cite{Munoz2012} extend the method for robust disturbance rejection and the simultaneous consideration of steady state stability and disturbance rejection, respectively. 
\cite{Kastsian2010} cover the case of fixed points of discrete time systems.
In the present paper, we extend the normal vector approach to stability constraints for periodic solutions of ODE systems. Similar but preliminary results are reported in \cite{Kastsian2012}.

We summarize some related methods for optimization of periodic processes in the remainder of this section. We comment on their ability to cope with uncertain model parameters and stability boundaries where appropriate.
The question whether periodic operation improves the system performance can be answered with the $\pi$-criterion \citep{Sterman1990, Parulekar1998}. Application of the $\pi$-criterion results in an optimal frequency for a sinusoidal input, but the criterion does not provide any information on the optimal amplitude and it does not apply to other input types. \cite{Avino2006} show that  in some situations it can even provide misleading results. The parameter continuation method described by \cite{Avino2006} gives a precise optimal point, but it is difficult to apply continuation methods for models with more than, say, two or three optimization or uncertain parameters. 

\cite{Mombaur2005a, Mombaur2005} and \cite{Mombaur2009} optimize periodic motions by solving two-level optimization problems. 
They optimize the economic objective function and minimize the spectral radius at the first and second level, respectively. The authors guarantee the resulting periodic orbits to be stable by minimizing the spectral radius and forcing all eigenvalues to have moduli strictly smaller than one. Parametric uncertainties in the underlying process models are not considered.

\cite{Burke2003} suggest minimizing the pseudo-spectral radius to guarantee robust stability. The pseudo-spectral radius measures 
the largest modulus of the eigenvalues of matrices which vary in an $\epsilon$-neighborhood of the reference matrix. The $\epsilon$-neighborhood is defined with the standard Euclidean norm. Since the pseudo-spectral radius typically is a nonsmooth function of the corresponding Jacobian entries, \cite{Vanbiervliet2009} and \cite{Diehl2009} proposed to use the smoothed spectral radius. The smoothed spectral radius is based on the $H_2$-norm and computed by solving  relaxed Lyapunov equations. When robustness is addressed with the pseudo-spectral radius or with the smoothed spectral radius it is difficult to consider parametric uncertainty.

\cite{Chang2011} consider parametric uncertainty for optimal steady state solutions and possible extension of the proposed method to oscillating processes. The authors solve semi-infinite programs, where stability constraints are addressed with the Routh-Hurwitz criterion. Note that the normal vector method proposed in the present paper does not use semi-infinite programs, but finite-dimensional nonlinear programs.  

The paper is organized as follows. We begin with a formal problem statement in Section \ref{sec:SystemClass} and outline of the normal vector method in Section \ref{sec:OutlineOfNVMethod}. In Section \ref{sec:CritMan} the characterization of the stability boundaries, or more generally critical manifolds, is introduced. 
The normal vectors to these critical manifolds and the nonlinear programs based on them are discussed in Section \ref{sec:NVOptim}. The proposed method is illustrated in Section \ref{sec:Examp}.  A conclusion is stated in Section \ref{sec:Concl}.

\section{System class and optimization problems of interest} \label{sec:SystemClass}
We consider dynamic systems described by a set of nonlinear parameterized ordinary differential equations 
\begin{eqnarray} \label{eq:sysODEs}
\dot{x}(t) = f(x(t),\alpha), \quad x(0)=x_{0},
\end{eqnarray}
where $x(t) \in \R^{n_{x}}$ and $\alpha \in 
\R^{n_{\alpha}}$ denote state variables and parameters, respectively. The function $f$ maps from some open subset of $\R^{n_{x}} \times \R^{n_{\alpha}}$ into $\R^{n_{x}}$ and is assumed to be smooth with respect to all variables and parameters. 

The simplest solutions of (\ref{eq:sysODEs}) are the equilibria, i.e., points $(x, \alpha)\in\R^{n_x}\times\R^{n_\alpha}$ such that  
\begin{eqnarray} \label{eq:SteadyState}
f(x,\alpha)=0.
\end{eqnarray}
The second class of solutions of ODE systems (\ref{eq:sysODEs}) that we consider are periodic orbits $(x(t),T,\alpha)$. Periodic orbits are solutions of (\ref{eq:sysODEs}) that satisfy the additional boundary condition
\begin{eqnarray} \label{eq:SysGenWithPeriod}
x(0)-x(T)=0,
\end{eqnarray} 
where the smallest admissible $T>0$ is the period of the orbit.
We are interested in finding equilibria or periodic solutions that are optimal with respect to a real valued objective function $\phi$, which may represent product concentration, productivity, or economic profit, for example. 
The optimal periodic solution  is determined by solving 
the optimization problem
\begin{eqnarray}
&{\max\limits_{{x^{(0)}(t)},T^{(0)}, \alpha^{(0)}}}  & {\phi(x^{(0)}(t),T^{(0)},\alpha^{(0)})} \nonumber \\
&{\rm s. t.} & \dot{x}^{(0)}(t) = f(x^{(0)}(t),\alpha^{(0)}),\nonumber\\ [-2.5 ex] \label{eq:OptUnconstrainedPer} \\ [-2.5ex]
&{} & {0} = {x^{(0)}(0)-x^{(0)}(T^{(0)})},\nonumber\\
&{} & {0} \leq {h(x^{(0)}(t),T^{(0)},\alpha^{(0)})}.\nonumber
\end{eqnarray}
The optimal equilibrium solution is found by solving
\begin{eqnarray}\label{eq:OptUnconstrainedSt}
&{\max\limits_{x^{(0)}, \alpha^{(0)}}}  & {\phi(x^{(0)},\alpha^{(0)})} \nonumber \\
&{\rm s. t.}  & {0} = {f(x^{(0)},\alpha^{(0)})},\\
&{} & {0} \leq {h(x^{(0)},\alpha^{(0)})}.\nonumber
\end{eqnarray}
We denote the objective function $\phi$ and the inequality constraints $h$ by the same symbols in both cases \eqref{eq:OptUnconstrainedPer} and \eqref{eq:OptUnconstrainedSt} for simplicity.           
The inequalities $h\ge 0$ model physical or economic constraints. Functions $\phi$ and  $h$  map from an open subset of  
$\mbox{$\R^{n_{x}}\times \R^{+}\times \R^{n_{\alpha}}$}$\,or $\mbox{$\mathbb{R}^{n_x}\times \mathbb{R}^{n_\alpha}$}$\,into $\R$ and $\R^{n_h}$,
 respectively, and are assumed to be smooth with respect to all variables and parameters. 
Note that we sometimes have to solve both optimization problems (\ref{eq:OptUnconstrainedPer}) and (\ref{eq:OptUnconstrainedSt}) and compare their objective function values to decide whether periodic or steady state operation is optimal.

\section{Outline of the normal vector approach}\label{sec:OutlineOfNVMethod}
\begin{figure}
\centering
\includegraphics{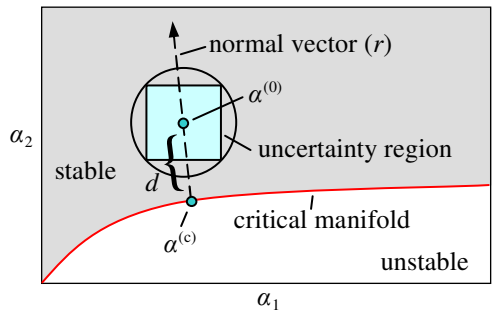}
\caption{Manifold of critical points for a hypothetical model. The parametric distance $d$ between $\alpha^{(0)}$ and the critical boundary can be measured along the vector $r$, which is normal to the manifold of critical points and passes through the candidate optimal point $\alpha^{(0)}$.} \label{fig:NVMethodSketch}
\end{figure}
The central idea of the normal vector method is that the parametric distance between the optimal point and a critical boundary can be measured along the normal direction to this boundary \citep{Dobson1993}. This idea is sketched in Figure \ref{fig:NVMethodSketch}. By ``critical boundary" we refer to boundaries in the space of the parameters $\alpha$ that separate regions with different dynamical properties of the system from one another. 
Typical critical boundaries separate regions with stable modes of operation from unstable ones. 
In this case the boundary is a projection of a manifold of bifurcation points onto the parameter space (see, e.g., \citet{Kuznetsov1998, Seydel1988}).

Figure \ref{fig:NVMethodSketch} illustrates how to force a candidate optimal point into the region with the desired dynamical properties. Essentially, the distance $d$ between the candidate optimal point $\alpha^{(0)}$ and the closest point $\alpha^{(c)}$ on the critical boundary must be sufficiently large  \citep{Moennigmann2002}.
This requirement can be enforced with the constraints 
\begin{eqnarray} \label{eq:GenNVConstraints}
\alpha^{(0)}-\alpha^{(c)}- d \frac{r}{\|r\|} = 0, \quad
d- d_{\min} \geq 0, 
\end{eqnarray}
where $\alpha^{(0)} \in \R^{n_\alpha}$ refers to the parameter values of the candidate optimal point, $\alpha^{(c)} \in \R^{n_\alpha}$ denotes the point on the critical boundary to which the normal vector $r \in \R^{n_\alpha}$ is stated, $d \in \R$ is the distance between $\alpha^{(0)}$ and the critical boundary, and $||\cdot||$ is the Euclidean norm.  
The choice of $d_{\min}$ will be explained below. If more than one critical manifold exist, or one or more critical manifolds are nonconvex, multiple constraints of the type (\ref{eq:GenNVConstraints}) have to be stated. This is detailed in Section \ref{sec:NVOptim}. We refer to constraints of the form~\eqref{eq:GenNVConstraints} as ''normal vector constraints``.  

It remains to take uncertain parameters in the model (\ref{eq:sysODEs}) into account.
We assume that the parameters $\alpha_i$ lie in intervals 
\begin{eqnarray}\label{eq:UncRegion}
\alpha_{i} \in [\alpha_{i}^{(0)}-\Delta \alpha_{i},\alpha_{i}^{(0)}+\Delta
\alpha_{i}], \quad i=1,\dots,n_\alpha,
\end{eqnarray}
where $\alpha_{i}^{(0)}$ are the central values of the independent uncertainty intervals and $\Delta \alpha_{i}$ represent the uncertainties. Since the parameters $\alpha_i$ may not have the same physical unit we introduce a simple metric to measure distances in the parameter space. Specifically, we measure the parameters in units of their  uncertainty $\Delta \alpha_{i}$. This is equivalent to  rescaling (\ref{eq:UncRegion}) according to 
\begin{eqnarray}\label{eq:ParScaling}
  \alpha_{i} \rightarrow \frac{\alpha_{i}}{\Delta
    \alpha_{i}},~\alpha^{(0)}_{i} \rightarrow \frac{\alpha^{(0)}_{i}}{\Delta
    \alpha_{i}}. 
\end{eqnarray}
The uncertainty region (\ref{eq:UncRegion}) then reads as
\begin{eqnarray}\label{eq:UncRegionScaled}
  \frac{\alpha_{i}}{\Delta \alpha_{i}} \in
  [\frac{\alpha^{(0)}_{i}}{\Delta \alpha_{i}}-1,\frac{\alpha^{(0)}_{i}}{\Delta
    \alpha_{i}}+1]\quad {\rm for~} i=1,\cdots,n_{\alpha}. 
\end{eqnarray}
In Figure \ref{fig:NVMethodSketch} and in what follows we assume that parameters $\alpha$ and $\alpha^{(0)}$ are scaled according to (\ref{eq:ParScaling}). The uncertainty region (\ref{eq:UncRegionScaled}) is sketched in Figure \ref{fig:NVMethodSketch}. It can be overestimated by a hyperball of radius $d_{\min}= \sqrt{n_{\alpha}}$.
The circle in Figure \ref{fig:NVMethodSketch} illustrates the two-dimensional case, i.e., $n_{\alpha}= 2$.  
Less conservative approximations for the uncertainty region (\ref{eq:UncRegionScaled}) than a hyperball exist but are not used here for simplicity. For details we refer the reader to \cite{Gerhard2008, Kastsian2010}. 

A second choice of the minimal distance $d_{\min}$ to the critical boundary is $d_{\min}=0$. 
In this case the candidate optimal point may lie on the critical boundary. 
If an uncertainty region (\ref{eq:UncRegionScaled}) is considered, some points of operation in the uncertainty region may cross the critical boundary.

\section{Critical manifolds of ODE systems with periodic solutions}\label{sec:CritMan}
Section \ref{sec:CritBoundsIllustration} reviews some notions from nonlinear systems theory (see, e.g., \citet{Kuznetsov1998} or \citet{Seydel1988}). The types of bifurcation points needed to describe the stability boundaries are summarized in Section \ref{subsec:StabilityBoundaries}. These boundaries are illustrated with a model of a peroxidase-oxidase reaction system in Section \ref{sec:ModelDescription}. 

\subsection{Stability analysis of periodic orbits}\label{sec:CritBoundsIllustration}
\begin{figure}
\centering
\subfloat[]{\includegraphics{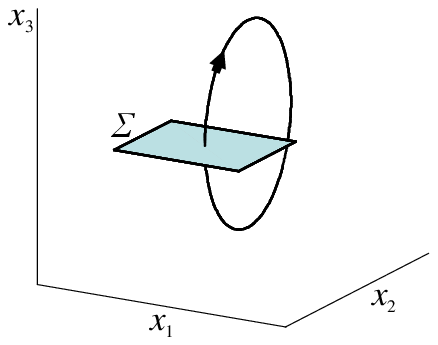} \label{fig:PoinacareMapSub1}}
\subfloat[]{\includegraphics{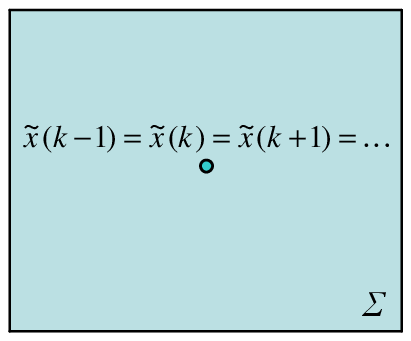} \label{fig:PoinacareMapSub2}}\\
\subfloat[]{\includegraphics{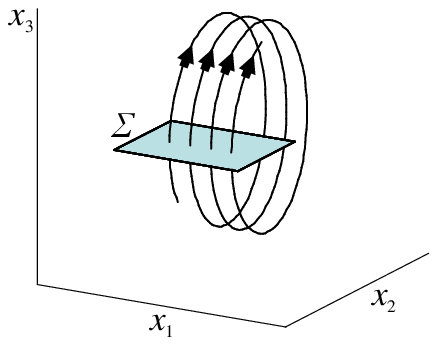} \label{fig:PoinacareMapSub3}}
\subfloat[]{\includegraphics{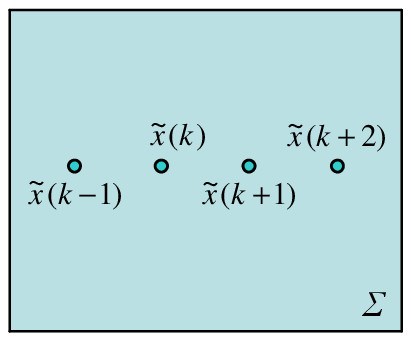} \label{fig:PoinacareMapSub4}}
\caption{Sketch of a periodic orbit (parts a, b) and a disturbed periodic orbit (parts c,d) of a three-dimensional system. In any transversal plane $\Sigma$ to the orbit, the periodic orbit appears as a fixed point of a discrete time system. Stability properties of the periodic orbit can conveniently be investigated by analyzing the stability properties of this discrete time system.} 
\label{fig:PoincareMap}
\end{figure}

  We briefly introduce the Poincar\'e section and Poincar\'e map, since they are instrumental for describing the stability properties of periodic orbits. See Figure \ref{fig:PoincareMap} for an illustration. The situation sketched in Figure \ref{fig:PoincareMap} can be described more specifically as follows. 
  
Let $\f(x_0,t,\alpha)$ denote the solution of (\ref{eq:sysODEs}) at time $t$ for the initial condition $x(0)=x_{0}$. Assume this solution is a periodic orbit with period $T$. It therefore satisfies Equation (\ref{eq:SysGenWithPeriod}), i.e., 
\begin{eqnarray}\label{eq:PerCondForPhi}
\f(x_0,T,\alpha) -x_0 = 0.
\end{eqnarray}
Since the Poincar\'e section $\Sigma$ shown in Figure \ref{fig:PoincareMap} can be shifted to intersect the orbit $\f$ at any other point, it is not unique. 
A particular Poincar\'e section is uniquely defined by specifying the point of its intersection with the periodic orbit $\f$, and requiring $\Sigma$ to be transversal (orthogonal) to the tangent to $\f$ at this point.
Formally, this is equivalent to introducing a phase condition
\begin{eqnarray}\label{eq:PhaseCondition}
s(x_0,T,\alpha)=0,
\end{eqnarray}
where $s$ maps from a subset of $\R^{n_x}\times \R^{+} \times \R^{n_\alpha}$ into $\R$. A discussion of the phase condition (\ref{eq:PhaseCondition}) is beyond the paper. We refer the reader to \cite{Kuznetsov1998} for details. Without restriction we choose the initial condition $x_0$ as the point of intersection. 

The Poincar\'e map $\Pi$ is the function that maps a point $x(kT) \in \Sigma$, $k= \{0, 1, 2, \dots\}$ onto the point $x((k+1)T)$ attained along the periodic orbit \eqref{eq:PerCondForPhi} after one period $T$, i.e., 
\begin{eqnarray}\label{eq:PoincareMap}
 \Pi: \Sigma \rightarrow \Sigma, \quad 
 x(kT) \rightarrow x((k+1)T)= \f(x(kT), T, \alpha).
\end{eqnarray}
By a slight abuse of notation we denote $x(kT)$, $x((k+1)T)$, etc.\ by $x(k)$, $x(k+1)$, respectively, to stress that the Poincar\'e map yields a discrete time system. The Poincar\'e map is usually defined in local coordinates $\tilde{x}=(\tilde{x}_1,\dots,\tilde{x}_{n_x-1})\in\R^{n_x- 1}$ on $\Sigma$ (see Figures \ref{fig:PoinacareMapSub2} and \ref{fig:PoinacareMapSub4}).  
This results in a discrete time system of the form
\begin{eqnarray}\label{eq:LocalPoincareMap}
  \Pi:\R^{n_x-1} \rightarrow \R^{n_x-1}, \quad
  \tilde{x}(k)\rightarrow \tilde{x}({k+1})=\Pi(\tilde{x}(k)),
\end{eqnarray}
where we use the same symbol $\Pi$ in \eqref{eq:PoincareMap} and \eqref{eq:LocalPoincareMap} for simplicity.

The orbit $\f$ is a periodic orbit of the continuous time system \eqref{eq:sysODEs} if and only if the intersection point $x_0$ is a fixed point of the Poincar\'e map, i.e. $\tilde{x}_0= \Pi(\tilde{x}_0)$, where $\tilde{x}_0$ denotes $x_0$ expressed in the local coordinates. 
The stability of the periodic orbit $\f$ can be investigated by analyzing the stability of the corresponding fixed point of the discrete time system.
More precisely, let $\f_{x_0}$ denote the Jacobian matrix of $\f(x_0, T, \alpha)$ with respect to $x_0$. 
This Jacobian evaluated at the fixed point is often referred to as the monodromy matrix. We denote it by
\begin{eqnarray}\label{eq:MonMatr}
M=\f_{x_0}(x_{0},T,\alpha),
\end{eqnarray}
for brevity. $M$ has eigenvalues $\lambda= 1$, $\lambda_1, \dots, \lambda_{n_x-1}$ if the Poincar\'e map \eqref{eq:LocalPoincareMap} has eigenvalues $\lambda_1, \dots, \lambda_{n_x- 1}$. The periodic solution $\f$ is locally asymptotically stable if all eigenvalues $\lambda_i$ are strictly inside the unit circle, or equivalently
\begin{eqnarray}\label{eq:EigenvalueCondition}
  |\lambda_i|<1 \mbox{ for all } i
\end{eqnarray}
(see, e.g., \cite{Kuznetsov1998}).

\subsection{Stability boundaries}\label{subsec:StabilityBoundaries}
  Consider the periodic orbit $\f(x_0, T, \alpha)$ of \eqref{eq:sysODEs} introduced in \eqref{eq:PerCondForPhi} again. This periodic orbit exists for certain fixed values of the parameter $\alpha$. If we change one or more of these parameters slightly, we expect the periodic orbit and the eigenvalues of the Poincar\'e map \eqref{eq:LocalPoincareMap} to vary slightly and continuously only. In particular we expect the eigenvalues to stay strictly inside the unit circle, and hence the periodic orbit to remain stable, for sufficiently small changes in $\alpha$. 
  If we intent to find an optimal periodic orbit, we generally have to admit large changes in $\alpha$, however. This implies that one or more eigenvalues $\lambda_i$ may leave the unit circle thus causing a loss of stability. 
  Bifurcation theory distinguishes three cases for such a loss of stability to occur, because each of these cases results in a particular change in the fixed point or periodic orbit behavior (see, e.g., \cite{Kuznetsov1998}). 
  At a Neimark-Sacker (torus) bifurcation point, a pair $(\lambda_{n_1},\lambda_{n_2})$ of complex conjugate eigenvalues of $M$ \eqref{eq:MonMatr} appears on the unit circle, $\lambda_{n_1}=e^{i \theta}$ and $\lambda_{n_2}=e^{-i \theta}$. Flip (period doubling) bifurcation points are associated with an eigenvalue of $M$ equal to $\lambda_{p_1}=-1$, whereas for fold (saddle-node) bifurcation points $\lambda_{l_1}=1$.   
We treat the stability boundaries associated with Neimark-Sacker and flip bifurcations of cycles. Fold bifurcations of cycles can be treated accordingly, but are not considered here, since they do not appear in the examples 
in the following sections.

For periodic operation we will accept only stable orbits.
In the case of Neimark-Sacker bifurcation a periodic orbit becomes unstable with appearance of a quasiperiodic motion. 
A quasiperiodic orbit has a periodic pattern but with irregular components. 
In contrast to stable periodic orbits, quasiperiodic orbits do not return to their initial conditions. Furthermore, a transition from quasiperiodic to chaotic behavior can occur.  In the case of flip bifurcation, a periodic orbit loses stability through period doubling and a chaotic behavior can result from a series of period doublings.

The stability properties of steady states \eqref{eq:SteadyState} can be characterized in a similar fashion with the eigenvalues of the Jacobian $f_x$. An equilibrium $(x^{(0)},\alpha^{(0)})$ is locally asymptotically stable, if the real parts of all eigenvalues of the Jacobian $f_x(x^{(0)},\alpha^{(0)})$ are negative. Generally, for equilibria there are two ways how stability can be lost while varying the system parameters. Hopf bifurcation points arise with appearance of two complex conjugate, purely imaginary eigenvalues. Saddle-node bifurcation points are associated with a real zero eigenvalue.

For steady state operation points we will also require stability. In the case of Hopf bifurcation the transition to either stable or undamped oscillations appears.
As we mentioned above stable oscillations will be permitted for the process operation. However, undamped oscillations are undesired. 
Saddle-node bifurcations lead to infeasible regions, where no solutions exist.
All cases of bifurcations of equlibria and cycles can be treated with normal vector constraints \eqref{eq:GenNVConstraints}, where $d_{\min}$ is chosen properly. 

We use the abbreviations ''NS``, ''flip``, ''sn``, and ''Hopf`` in figures and equations to refer to Neimark-Sacker, flip, saddle-node, and Hopf bifurcation points, respectively. Saddle-node bifurcations exist for both periodic orbits and equilibria. Here ''Saddle-node`` and ''sn`` always refer to the equilibrium case if not noted otherwise.

\subsection{The peroxidase-oxidase reaction model}\label{sec:ModelDescription}
We introduce the peroxidase-oxidase reaction model that will later be optimized in Section \ref{sec:NADExamp}. 
The model is introduced here already, because it can be used in illustrations throughout the paper this way. 

The peroxidase-oxidase reaction model describes the aerobic oxidation of nicotinamide adenine dinucleotide hydrid (NADH) by molecular oxygen, which is catalyzed by horseradish peroxidase enzyme (HRP). The overall net reaction is given by
\begin{eqnarray}\label{eq:NADHreaction}
{\rm 2NADH + O_2 + 2H^{+} \xrightarrow{HRP} 2NAD^{+} + 2H_2O.}
\end{eqnarray}
The reaction takes place in the presence of methylene blue and 2,4-dichloro\-phenol \citep{ Steinmetz1993, Larter2003}. 
The peroxidase-oxidase reaction plays an important role in the production of lignin, a polymer that makes wood hard \citep{Halliwell1978, Maeder1982}. 
The reaction product $\rm{NAD^{+}}$ is also of interest in pharmacology \citep{Khan2007,Sauve2008}.

There exists no universally agreed mathematical model for the peroxidase-oxidase reaction (\ref{eq:NADHreaction}), but the characteristics of this reaction have been effectively modeled by using a simplified eight-step mechanism proposed by \cite{Olsen1993}
\begin{subequations} \label{eq:ChemReactNADModel}
\begin{eqnarray}\nonumber
\begin{tabular}{rcrcrcr}
$B+X \xrightarrow{k_1} 2X,$ & & (\ref{eq:ChemReactNADModel}a) & \quad\quad\quad & $Y \xrightarrow{k_5} Q,$ &  & (\ref{eq:ChemReactNADModel}e)\\
$2X \xrightarrow{k_2} 2Y,\,$ & & (\ref{eq:ChemReactNADModel}b) & & $X_0 \xrightarrow{k_6} X,$ & & (\ref{eq:ChemReactNADModel}f)\\
$A+B+Y \xrightarrow{k_3} 3X,$ & & (\ref{eq:ChemReactNADModel}c) & & $A_0 \begin{array}{c}  \xrightarrow{k_7} \\ 
[-2.5 ex] \xleftarrow[k_7]{} \end{array} A,$ & & (\ref{eq:ChemReactNADModel}g)\\
$X  \xrightarrow{k_4} P,~\,$ & & (\ref{eq:ChemReactNADModel}d) & & $B_0 \xrightarrow{k_8} B.$ & & (\ref{eq:ChemReactNADModel}h)
\end{tabular} 
\end{eqnarray}
\end{subequations}
$A$ and $B$ denote the concentrations of the reactants O$_2$ and NADH, respectively. $A_0$ and $B_0$ are the concentrations of $A$ and $B$ in the feed streams, respectively. $P$ and $Q$ are the reaction products. 
$X$ and $Y$ represent intermediate free radicals NAD$^{\bullet}$
and oxyferrous peroxidase \citep{Aguda1989}, respectively. Note that NAD$^{\bullet}$ denotes electrically neutral radicals of nicotinamide adenine dinucleotide and oxyferrous peroxidase is sometimes called compound ${\rm \MakeUppercase{\romannumeral 3}}$ \citep{Aguda1989}. 

The steps (\ref{eq:ChemReactNADModel}a) and (\ref{eq:ChemReactNADModel}b)--(\ref{eq:ChemReactNADModel}c) form two routes for the autocatalytic production of intermediate NAD$^{\bullet}$. Reaction (\ref{eq:ChemReactNADModel}d) and (\ref{eq:ChemReactNADModel}e) are two linear radical termination steps, while reaction (\ref{eq:ChemReactNADModel}f) is the initialization step of the radicals.  The equilibrium between gaseous O$_{2}$ and the liquid phase is addressed in (\ref{eq:ChemReactNADModel}g). Step (\ref{eq:ChemReactNADModel}h) refers to the inflow of NADH.

The following model results from applying the law of mass action to the reaction mechanism~\eqref{eq:ChemReactNADModel} 
\begin{eqnarray}
\dot{A} &=& k_7(A_0-A)-k_3 ABY, \nonumber\\
\dot{B} &=& k_8 B_0 -k_1 BX - k_3 ABY, \nonumber\\ [-2.5 ex] \label{eq:SysNAD} \\ [-2.5ex]
\dot{X}&=& k_1 BX- 2k_2 X^2 +3k_3 ABY -k_4 X+ k_6 X_0, \nonumber\\ 
\dot{Y}&=& 2k_2 X^2 -k_3 ABY -k_5 Y,\nonumber
\end{eqnarray}
where all variables are dimensionless \citep{Olsen1993}.
  The parameters $k_1$ and $k_3$ define the total peroxidase enzyme concentration and the concentration of 2,4-dichlorophenol \citep{Steinmetz1993, Larter2003}. 
  When maximizing the NAD$^{+}$ concentration in Section \ref{sec:NADExamp}, we search for the optimal values for these parameters $k_1$ and $k_3$ within the bounds
  \begin{eqnarray}\label{eq:BoundsNADproduction}
  0.1 \leq k_1 \leq 0.5, \quad 0.001 \leq k_3 \leq 0.05 
  \end{eqnarray}
  Following \cite{Steinmetz1993} and \cite{Larter2003}, where the model (\ref{eq:SysNAD}) is verified in laboratory experiments, we assume the exact optimal values for $k_1$ and $k_3$ cannot be controlled to arbitrary precision, but they may drift within certain error bounds. Consequently, $k_1$ and $k_2$ are uncertain parameters. The uncertainty is stated precisely in \eqref{eq:UncRegOlsen} below. 
  The other parameters are fixed to the values $k_2=250$, $k_4=20$, $k_5=5.35$, $k_6 X_0=10^{-5}$, $k_7=0.1$, $k_8 B_0=0.825$, and $A_0=8$ \citep{Steinmetz1993}.

\begin{figure}
\centering
\includegraphics{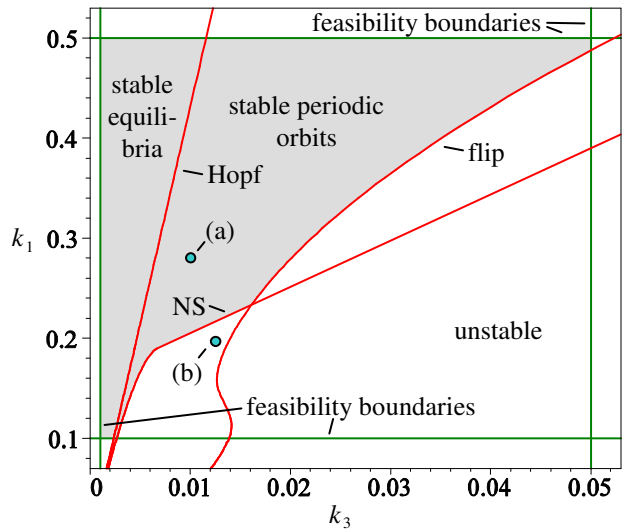}
\caption{Critical boundaries for the peroxidase-oxidase reaction model.}
\label{fig:BifPicPOmodel}
\end{figure}
\begin{figure}
\centering
\subfloat[]{\includegraphics{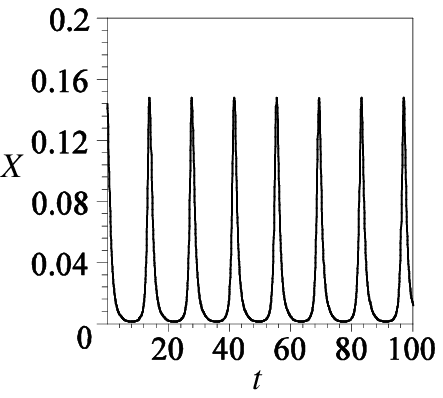} \label{fig:TimeVsXSubA}}
\subfloat[]{\includegraphics{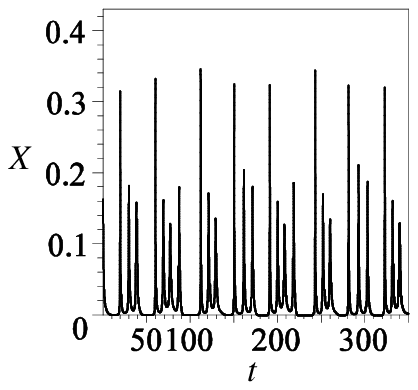} \label{fig:TimeVsXSubB}}\\
\subfloat[]{\includegraphics{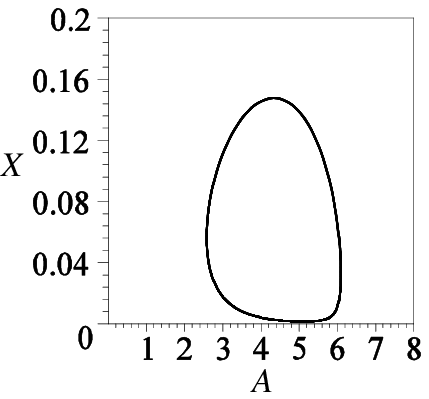} \label{fig:TimeVsXSubC}}
\subfloat[]{\includegraphics{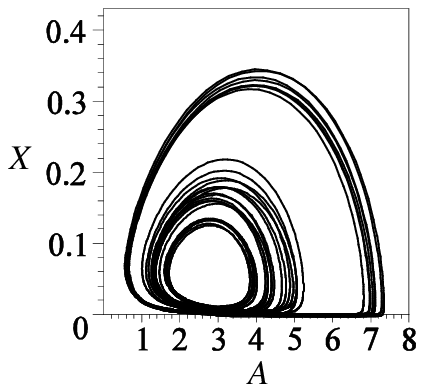} \label{fig:TimeVsXSubD}}
\caption{Time series and phase portraits from desired (diagrams a, c) and undesired (diagrams b, d) regions in Figure \ref{fig:BifPicPOmodel}}.
\label{fig:TimeVsX}
\end{figure}

Figure \ref{fig:BifPicPOmodel} shows the bifurcation points of the peroxidase-oxidase reaction model (\ref{eq:SysNAD}) in the plane spanned by the two uncertain parameters $k_1$ and $k_3$. 
The Hopf bifurcation points give rise to a stable periodic solution and an unstable equilibrium in this particular reaction system. 
The resulting stable periodic solutions lose stability at the Neimark-Sacker or flip bifurcation points of cycles. 
The lines labeled ``feasibility boundaries" result from the constraints \eqref{eq:BoundsNADproduction}. 
The regions in which the desired dynamical properties exist, i.e., stable and feasible equilibria or stable and feasible periodic orbits, are shaded in Figure~\ref{fig:BifPicPOmodel}. 
Diagrams \ref{fig:TimeVsXSubA} and \ref{fig:TimeVsXSubB} show time series evaluated at points labeled (a) and (b) in Figure~\ref{fig:BifPicPOmodel}, respectively. Diagrams \ref{fig:TimeVsXSubC} and \ref{fig:TimeVsXSubD} show the respective phase portraits.

\subsection{Systems of equations for critical manifolds} \label{sec:MathDescriptionOfBounds}
Critical boundaries like those shown in Figure \ref{fig:BifPicPOmodel} can be
described with so-called augmented systems. Essentially, augmented systems formally state the necessary conditions for the critical eigenvalues explained in Section \ref{subsec:StabilityBoundaries}.
The augmented systems are the basis for the calculation of the normal direction $r$ introduced in Figure \ref{fig:NVMethodSketch} and Equation~(\ref{eq:GenNVConstraints}).

We briefly explain the augmented system for flip bifurcation points of periodic solutions. This system reads as the following set of $2n_x+2$ equations \citep{Lust1997, Lust1999, Khinast2000}:
\begin{eqnarray}\label{eq:FlipSystem}
M^{(\rm{flip})}(x_0,T,\alpha):=\left(\begin{array}{c}
\f(x_0,T,\alpha) -x_0\\
s(x_0,T,\alpha)\\
\f_{x_0}(x_0,T,\alpha)w + w\\
w^{T}w - 1
\end{array}\right)=0,
\end{eqnarray}
where the first two lines are 
the periodicity and phase conditions discussed in Section \ref{subsec:StabilityBoundaries}, the third line ensures that the Poincar\'e map and the monodromy matrix \eqref{eq:MonMatr} have an eigenvalue $-1$ with eigenvector $w \in \R^{n_x}$, and the last line is the normalization of $w$. The system (\ref{eq:FlipSystem}) is nonsingular with respect to $x_0$, $T$, $w$, and one component of $\alpha \in \R^{n_\alpha}$, say $\alpha_1$, at any nondegenerate flip bifurcation point of cycles \citep{Lust1997, Lust1999}. 

The augmented system for Neimark-Sacker bifurcation points of periodic solutions is stated in \ref{app:NSAugSys} for completeness.
Augmented systems for saddle-node and Hopf bifurcations of equilibria of ODEs are omitted for brevity (see, e.g., \cite{Kuznetsov1998,Mangold2000,Beyn2002,Moennigmann2002}).

\section{Optimization with guaranteed robust stability} \label{sec:NVOptim}
  It remains to incorporate the critical boundaries described in the previous section into the process optimization problems (\ref{eq:OptUnconstrainedPer}) and (\ref{eq:OptUnconstrainedSt}). This is done with constraints of the form \eqref{eq:GenNVConstraints} for which the normal vector $r$ sketched in Figure \ref{fig:NVMethodSketch} is instrumental.
  We describe the systems of equations that define $r$ in Section \ref{sec:NVDirect}. We give only a brief description, since these equations are technical. Subsequently, we state the optimization problems with normal vector constraints for robust stability in Section \ref{sec:NVProbAug}.

\subsection{Systems of equations for the normal vectors}\label{sec:NVDirect}
Normal vector systems can be derived by applying the scheme of derivation proposed in \cite{Moennigmann2002} to the augmented systems. The details for the derivation of the particular systems treated here can be found in \cite{KastsianPhDThesis}.
We state the normal vector system for the case of the flip bifurcation of cycles as an example:
\begin{eqnarray}\nonumber
G^{(\rm{flip})}(p,\bar{x}^{(\rm{flip})},r):=
\left(\begin{array}{c}
M^{(\rm{flip})}(p)\\
\f_{x_0}^{T}(p) v + v +\gamma_1 w\\
v^{T}w-1\\
\f_{x_0}^{T}(p) u - u + s_{x_0}^{T}(p) \varkappa +  v^{T} \f_{x_0x_0}(p) w\\
\f_{T}^{T}(p) u + s_T(p) \varkappa +  v^{T} \f_{x_0T}(p) w\\
r - \f_{\alpha}^{T}(p) u - s_{\alpha}^{T}(p)  \varkappa -  v^{T} \f_{x_0 \alpha}(p) w
\end{array}\right)=0,
\end{eqnarray}
where
\begin{eqnarray}\label{eq:Abbreviations}
p=(x_0,T,\alpha),\quad\bar{x}^{(\rm{flip})}= (w, v, u, \varkappa, \gamma_1)
\end{eqnarray}
are introduced for brevity.
$M^{(\rm{flip})}(p)$ refers to Equation \eqref{eq:FlipSystem}.
The matrices $\f_{x_0} \in \R^{n_x \times n_x}$,  $s_{x_0} \in \R^{1 \times n_x}$, $\f_{\alpha} \in \R^{n_x \times n_\alpha}$,  and $s_{\alpha} \in \R^{1 \times n_\alpha}$ are the obvious matrices of derivatives with respect to $x_0$ and $\alpha$, respectively. Furthermore, $\f_T=f(x_0,\alpha) \in \R^{n_x}$, $\f_{x_0T}=f_x(x_0,\alpha) \in \R^{n_x \times n_x}$ and $s_T \in \R$ denote derivatives with respect to period $T$. Note that $\f_{x_0x_0}$ and $\f_{x_0\alpha}$ represent second order derivatives. 
The variables $w \in \R^{n_x}$ and $v \in \R^{n_x}$ are the eigenvectors of $\f_{x_0}$ and  its transpose $\f_{x_0}^T$, respectively, that correspond to eigenvalue $-1$. The symbol $r \in \R^{n_\alpha}$ denotes the normal vector.
Finally, $u \in \R^{n_x}$, $\varkappa \in \R$, and  $\gamma_1 \in \R$ are auxiliary variables.

The derivatives $\f_{x_0}$, $\f_{\alpha}$,  $\f_{x_0x_0}$, and $\f_{x_0 \alpha}$ can be obtained  with automatic differentiation. We use TIDES \citep{Abad2009} to calculate partial derivatives of $\f$. The other derivatives $\f_T$, $\f_{x_0T}$, $s_{x_0}$, $s_{T}$, and $s_\alpha$  
can be obtained with symbolic differentiation. 
The normal vector systems for the remaining bifurcations points are stated in the appendix. 

In general, the normal vector systems for manifolds of bifurcation points of cycles have the form
\begin{eqnarray}\label{eq:CritManForCycles}
G^{(c)}(p,\bar{x}^{(c)},r)=0,
\end{eqnarray}
where $c\in\{\rm{NS, flip}\}$ indicates the type of the bifurcation and normal vector system.
The general form of the normal vector systems for bifurcations of equilibria is 
\begin{eqnarray}\label{eq:CritManForEquil}
G^{(c)}(q,\bar{x}^{(c)},r)=0,
\end{eqnarray}
  where 
\begin{eqnarray}\label{eq:AbbreviationsII}
   q=(x,\alpha)
 \end{eqnarray}
is introduced for brevity and $c\in\{\rm{Hopf}, \rm{sn}\}$. 
Other types of critical manifolds (e.g., feasibility constraints) can be considered in the same manner, but are not necessary here.

\subsection{Optimization procedures with the normal vector constraints} \label{sec:NVProbAug}
  The region of stable behavior is generally bounded by more than one critical boundary. For example, Figure \ref{fig:BifPicPOmodel} shows there exist a Hopf, a flip, and a Neimark-Sacker boundary for the peroxidase-oxidase reaction model \eqref{eq:ChemReactNADModel}. 
  We assume there exist $i_{\max}$ critical boundaries and introduce tuples $(c_i, i)$, $c_i\in\{\rm{NS}, \rm{flip}, \rm{Hopf}, \rm{sn}\}$ to indicate the type of the respective boundary. 
  Without restriction we assume that the critical boundaries $1, \dots, \tilde{i}_{\max}$ and $\tilde{i}_{\max}+1, \dots, i_{\max}$ belong to periodic orbits and equilibria, respectively. 
Combining the optimization problem \eqref{eq:OptUnconstrainedPer} for periodic operation with the normal vector constraints \eqref{eq:CritManForCycles} and the defining system for $r$ from \eqref{eq:GenNVConstraints} results in the robust optimization problem
\begin{subequations}\label{eq:OptProbelmWithNVConstr}
\begin{eqnarray}
&{\max\limits_{{x^{(0)}(t)},T^{(0)}, \alpha^{(0)}}}  & {\phi(x^{(0)}(t),T^{(0)},\alpha^{(0)})} \nonumber \\
&{\rm s. t.}  & \dot{x}^{(0)}(t) = f(x^{(0)}(t),\alpha^{(0)}),\label{eq:ConstrNVForPer1}\\
&{} & 0=x^{(0)}(0)-x^{(0)}(T^{(0)}),\label{eq:ConstrNVForPer2}\\
&{} & 0 \leq h(x^{(0)}(t),T^{(0)},\alpha^{(0)}),\label{eq:ConstrNVPhys}\\
&{} & 0 = G^{(c_i,i)}(p^{(i)},\bar{x}^{(c_i,i)},r^{(i)}),\quad 
      i= 1, \dots, \tilde{i}_{\max}, \label{eq:FirstNVConstr}\\ 
&{} & 0 = G^{(c_j,j)}(q^{(j)},\bar{x}^{(c_j,j)},r^{(j)}),\quad 
      j= \tilde{i}_{\max}+1, \dots, i_{\max}, \label{eq:NVConstrEq1}\\  
&{} & 0 = \alpha^{(0)}-\alpha^{(c,k)}- d^{(k)} \frac{r^{(k)}}{\|r^{(k)}\|}, \quad 
      k= 1, \dots, i_{\max},
      \label{eq:NVConstrEq2}\\ 
&{} & {0} \leq {d^{(k)}}-d_{\min}^{(k)}, \quad 
      k= 1, \dots, i_{\max}. 
      \label{eq:LastNVConstr}
\end{eqnarray} 
\end{subequations} 
Constraints (\ref{eq:ConstrNVForPer1}) and (\ref{eq:ConstrNVForPer2}) ensure that the optimal solution  corresponds to a periodic orbit of the ODE system (\ref{eq:sysODEs}). 
Constraints (\ref{eq:ConstrNVPhys}) are the feasibility constraints from \eqref{eq:OptUnconstrainedPer}. 
Equations (\ref{eq:FirstNVConstr}) and (\ref{eq:NVConstrEq1}) state the normal vector systems (\ref{eq:CritManForCycles}) and (\ref{eq:CritManForEquil}), respectively.   
The symbol $r^{(k)}$ denotes the $k$th normal vector at point $\alpha^{(c_k,k)}$, which belongs to the $k$th critical boundary.
Constraints (\ref{eq:NVConstrEq2}) and (\ref{eq:LastNVConstr}) implement \eqref{eq:GenNVConstraints} for the $k$th critical boundary.

The corresponding augmented optimization problem for the optimal equilibrium 
$(x^{(0)},\alpha^{(0)})$ reads as
\begin{eqnarray}
&{\max\limits_{x^{(0)},\alpha^{(0)}}}  & {\phi(x^{(0)},\alpha^{(0)})} \nonumber \\
&{\rm s. t.}  & 0=f(x^{(0)},\alpha^{(0)}),\nonumber\\ [-2.5 ex] \label{eq:OptProbelmWithNVConstrEquil} \\ [-2.5ex] 
&{} & {0} \leq {h(x^{(0)},\alpha^{(0)})},\nonumber\\
&{} & \mbox{constraints (\ref{eq:FirstNVConstr})--(\ref{eq:LastNVConstr})},\nonumber
\end{eqnarray} 
where the first and second constraints are as in \eqref{eq:OptUnconstrainedSt}, and the normal vector constraints are adopted from (\ref{eq:OptProbelmWithNVConstr}).

If both equilibria and periodic orbits exist, we solve optimization problems (\ref{eq:OptProbelmWithNVConstr}) and (\ref{eq:OptProbelmWithNVConstrEquil}) and choose the maximum from the two resulting optimal modes of operation. 
The critical boundaries need not be known a priori, but can be automatically detected \citep{Moennigmann2007}.

\begin{figure}
\centering
\includegraphics{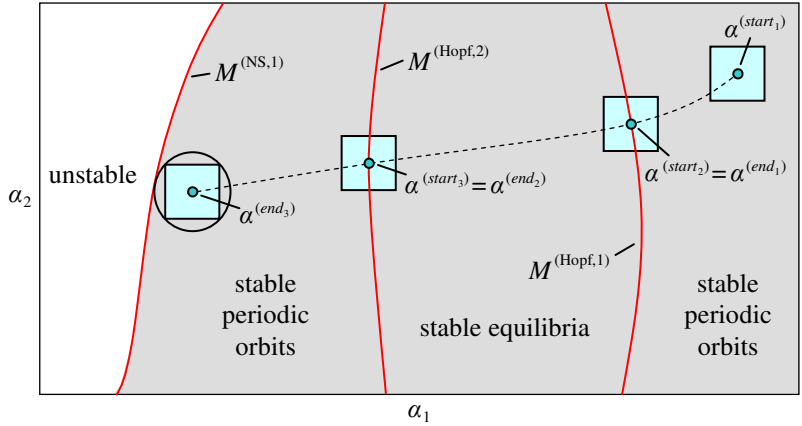}
\caption{Sketch of a situation in which switching between optimization problems (\ref{eq:OptProbelmWithNVConstr}) and  (\ref{eq:OptProbelmWithNVConstrEquil}) is necessary.}
\label{fig:BoundSwitch}
\end{figure}

If the optimal point from (\ref{eq:OptProbelmWithNVConstr}) or (\ref{eq:OptProbelmWithNVConstrEquil}) lies on a critical boundary that separates a region with periodic from a region with equilibrium solutions, i.e., on a Hopf bifurcation boundary, we have to carry out the optimization in both regions. Switching from one region to the other involves switching between problems (\ref{eq:OptProbelmWithNVConstr}) and (\ref{eq:OptProbelmWithNVConstrEquil}). Situations of this type are illustrated in Figure \ref{fig:BoundSwitch}. 
Assume we initialize the optimization problem (\ref{eq:OptProbelmWithNVConstr}) with a stable periodic orbit, which corresponds to $\alpha^{(\rm{start}_1)}$. 
By construction this optimization problem cannot cross a Hopf bifurcation boundary  
(manifold labeled $M^{(\rm{Hopf},1)}$ in Figure \ref{fig:BoundSwitch}). 
If such a boundary is encountered, 
an equilibrium that exist in the neighboring parameter region can be used to initialize 
optimization problem~(\ref{eq:OptProbelmWithNVConstrEquil}) ($\alpha^{(\rm{start}_2)}= \alpha^{(\rm{end}_1)}$ in Figure \ref{fig:BoundSwitch}). Conversely, the equilibrium optimization problem  (\ref{eq:OptProbelmWithNVConstrEquil}) cannot cross critical boundaries at which the equilibrium solution vanishes or becomes unstable (manifold labeled $M^{(\rm{Hopf},2)}$ in Figure \ref{fig:BoundSwitch}). Just as in the first case, a bifurcation point to a region with stable periodic behavior can be used to initialize a new optimization problem of the form (\ref{eq:OptProbelmWithNVConstr}) ($\alpha^{(\rm{start}_3)}=\alpha^{(\rm{end}_2)}$ in Figure \ref{fig:BoundSwitch}). Finally, there exist boundaries  (e.g., manifold labeled $M^{(\rm{NS},1)}$ in Figure \ref{fig:BoundSwitch}), at which no switching is required, since a stable equilibrium or periodic solution to the dynamical system (\ref{eq:sysODEs}) exists only on one side of the critical manifold.

For the solution of problems (\ref{eq:OptProbelmWithNVConstr}) and (\ref{eq:OptProbelmWithNVConstrEquil}) we use the SQP-solver NPSOL \citep{Gill2001} combined with the implicit Runge-Kutta method realization TWPBVPC \citep{Cash2005}. The gradient-based solver NPSOL requires the derivatives of the constraints of the optimization problems with respect to all optimization variables. 
 Therefore, third-order derivatives of $\f$ and $f$ are required. Analogously to the second-order derivatives, they can be determined for $\f$ with the automatic differentiation software TIDES \citep{Abad2009} and for $f$ with symbolic differentiation. Alternatively, the finite difference option of NPSOL \citep{Gill2001} can be used.

\section{Applications}\label{sec:Examp}
We apply the proposed method to two chemical reaction systems. Both systems exhibit autonomous oscillations and permit periodic or steady state operation. We note that switching between the periodic and equilibrium optimization problems \eqref{eq:OptProbelmWithNVConstr} and  \eqref{eq:OptProbelmWithNVConstrEquil} is necessary in the first application but not in the second one.
\subsection{Peroxidase-oxidase reaction model}\label{sec:NADExamp}
We optimize the peroxidase-oxidase reaction model (\ref{eq:SysNAD}) by maximizing the concentration of NAD$^{+}$. The objective function for equilibria reads as $\phi=X$. If the solution of (\ref{eq:SysNAD}) is a periodic orbit we maximize the average concentration
\begin{eqnarray}\label{eq:AverageConcentr}
\phi=\frac{1}{T}\int_{0}^{T}X(t) \mathrm dt.
\end{eqnarray}
For reference we first optimize the peroxidase-oxidase reaction model in Section \ref{sec:ResOptWithoutNV} without any stability constraints, i.e., we solve optimization problems (\ref{eq:OptUnconstrainedPer}) and (\ref{eq:OptUnconstrainedSt}) without normal vector constraints. The normal vector method is used in Section \ref{sec:ResOptWithNV}.

\subsubsection{Reference results obtained without normal vector constraints}\label{sec:ResOptWithoutNV}
\begin{figure}
\centering
\includegraphics{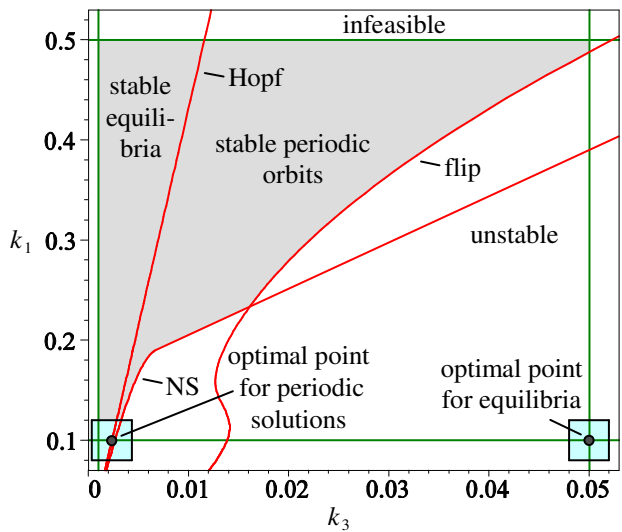}
\caption{Optimal robust points for the peroxidase-oxidase reaction model obtained without normal vector constraints.}
\label{fig:ResultsWithoutNV}
\end{figure}
Figure \ref{fig:ResultsWithoutNV} shows the optimal points that result from both the optimization problem (\ref{eq:OptUnconstrainedSt}) for equilibria, and from the optimization problem (\ref{eq:OptUnconstrainedPer}) for  periodic solutions of model (\ref{eq:SysNAD}). 
The uncertainty region (\ref{eq:UncRegion}) corresponds to 
\begin{eqnarray}\label{eq:UncRegOlsen}
  (k_1,k_3) \in ([k_1^{(0)}-\Delta k_1,k_1^{(0)}+\Delta k_1],[k_3^{(0)}-\Delta k_3,k_3^{(0)}+\Delta k_3]),
\end{eqnarray}
 where $\Delta k_1= 0.02$ and $\Delta k_3=0.002$. Symbols $k_1^{(0)}$ and $k_3^{(0)}$ denote optimization variables. 
The optimization results in an optimal but unstable equilibrium with $(k_1^{(0)},k_3^{(0)})=(0.1,0.05)$ 
and an objective function value $\phi=47.69 \cdot 10^{-3}$. Figure \ref{fig:ResultsWithoutNV} shows that the equilibrium is unstable, the entire robustness region (\ref{eq:UncRegOlsen}) lies in an unstable region, and a large fraction of it violates the boundaries $h\ge 0$.

Solving \eqref{eq:OptUnconstrainedPer} results in an optimal periodic orbit with 
an objective function value $\phi= 32.85 \cdot 10^{-3}$. Figure \ref{fig:ResultsWithoutNV} shows that the optimal periodic solution is stable but not robust, since a large fraction of the robustness region (\ref{eq:UncRegOlsen}) violates stability and feasibility boundaries.  

In summary, an optimization without stability constraints does not provide useful results for this sample process model. Both the optimal equilibrium and the optimal periodic orbit obtained from solving the optimization problems (\ref{eq:OptUnconstrainedPer}) and (\ref{eq:OptUnconstrainedSt}) are unacceptable from an operational point of view.

\subsubsection{Results of the robust optimization with normal vector constraints}\label{sec:ResOptWithNV}
In order to find the optimal stable and robust mode of operation, we force the optimal point to lie in the region where  stable equilibria or stable periodic orbits exist with the normal vector method. 
We start the optimization procedure  (\ref{eq:OptProbelmWithNVConstr}) with a stable periodic orbit. The minimal distances to critical boundaries of flip and Neimark-Sacker bifurcation points of cycles are set to $d_{\min}=\sqrt{n_\alpha}=\sqrt{2}$.  
The minimal distance to Hopf bifurcation points of equilibria is set to zero in order to permit switching from the optimization of periodic orbits to the optimization of equilibria or vice-versa. In fact, the optimization (\ref{eq:OptProbelmWithNVConstr}) drives the optimal point to the Hopf bifurcation boundary.
\begin{figure}
\centering
\includegraphics{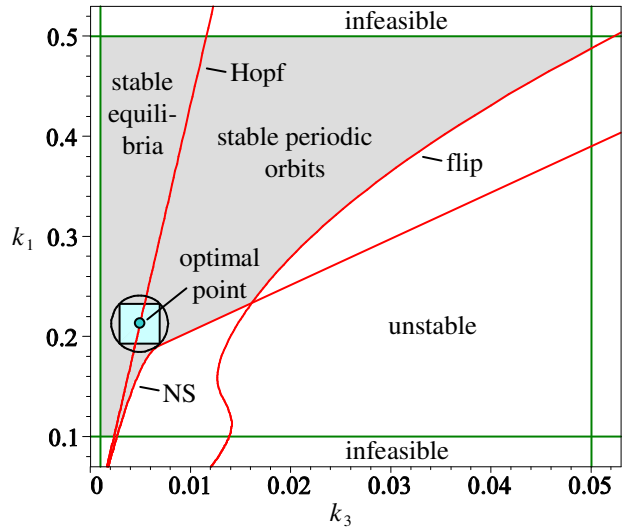}
\caption{Optimal robust point for the peroxidase-oxidase reaction model obtained with normal vector constraints.}
\label{fig:ResultsWithNVParam}
\end{figure}
Consequently, we switch to the problem (\ref{eq:OptProbelmWithNVConstrEquil}) that seeks for stable equilibria with the same values for $d_{\min}$.
The resulting optimal point, which is located on the Hopf boundary, is illustrated in Figure \ref{fig:ResultsWithNVParam}. 
The optimal parameter values are $(k_1^{(0)},k_3^{(0)})=(0.2126,0.00495)$. The objective function evaluates to $\phi=32.81 \cdot 10^{-3}$ at the optimal point.  
It is apparent from Figure \ref{fig:ResultsWithNVParam} that the entire robustness region around the optimal point lies in the stable region.  
The value of the objective function obtained is lower than those found in Section~\ref{sec:ResOptWithoutNV}, but we achieved stable and robust operation. 

\begin{figure}
\centering
\includegraphics{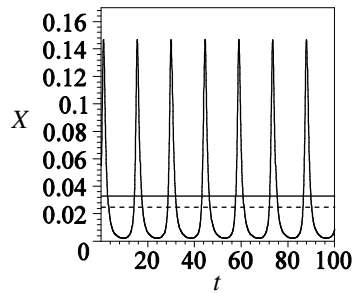}
\caption{Optimal robust steady state (solid line) and time series for the periodic solution for the closest point to the stability boundary ($(k_1, k_3)= (k_1^{(0)}-\Delta k_1,k_3^{(0)}+\Delta k_3))$. The dashed line indicates the average concentration of $X$ (\ref{eq:AverageConcentr}) for the periodic orbit.}
\label{fig:ResultsWithNVOrbit}
\end{figure}
The optimal equilibrium is depicted as a solid line in Figure \ref{fig:ResultsWithNVOrbit}.  For comparison we also show the periodic solution that results for parameter values $(k_1^{(0)}-\Delta k_1,k_3^{(0)}+\Delta k_3)$, i.e.,  the lower right corner of the robustness region. These values are chosen, because they correspond to the closest point to the stability boundary in the robustness region. The average concentration of $X$ for the periodic orbit is shown in Figure~\ref{fig:ResultsWithNVOrbit} for comparison. 

\subsection{Nonisothermal chemical reactor}
We consider an autocatalytic reaction  $P \rightarrow A  \rightarrow B  \rightarrow C$, where a relatively stable reactant $P$ is converted to a final product $C$ through two intermediate products $A$ and $B$ with the reaction steps
\begin{eqnarray}\label{eq:reactSchema}
\begin{tabular}{ r@{~$\rightarrow$~}l  l }
$P$ & $A$ & $\mbox{rate}$ $=k_0 p$, \\
$A$ & $B$ & $\mbox{rate}$ $=k_3 a$, \\
$A + 2B$ & $3B$ & $\mbox{rate}$ $=k_1 a b^2$,\\
$B$& $C + \rm{Heat}$ & $\mbox{rate}$ $=k_2 b.$
\end{tabular}
\end{eqnarray}
In (\ref{eq:reactSchema}) $k_0$ [s$^{-1}$], $k_1$ [m$^{6}\,$mol$^{-2}\,$s$^{-1}$], $k_2$ [s$^{-1}$], and $k_3$ [s$^{-1}$] are rate constants for the corresponding reactions. The symbols $p$, $a$, and $b$ denote the concentrations of $P$, $A$, and $B$, respectively, measured in [mol$\,$m$^{-3}$]. The model of (\ref{eq:reactSchema}) is adopted from \cite{Scott1990}.
The reaction rate equations for the concentrations $p$, $a$, and $b$, and the energy balance read as
\begin{eqnarray}
\dot{p} &=&-k_0 p, \nonumber \\  
\dot{a} &=& k_0 p -k_1 a b^2 -k_3 a, \nonumber \\[-2.5 ex]  \label{eq:ChemReactAB}  \\ [-2.5ex] 
\dot{b} &=& k_1 a b^2 +k_3 a - k_2 b, \nonumber\\
\dot{\tau} &=& \frac{1}{V c_p c_0}(V Q k_2 b - \chi S (\tau-\tau_a)),\nonumber
\end{eqnarray}
where $\tau$ [K] refers to the temperature, $V$ [m$^3$] denotes the reactor volume, $c_p$~[J$\,$mol$^{-1}\,$K$^{-1}$] is the molar heat capacity , $c_0$ [mol$\,$m$^{-3}$] is the molar density, $\chi$ [W$\,$m$^{-2}\,$K$^{-1}$)] is the surface heat transfer coefficient, $S$ [m$^2$] is the surface area, $Q$ [J$\,$mol$^{-1}$] is the heat of the exothermic reaction from (\ref{eq:reactSchema}), and $\tau_a$ [K] denotes the temperature of the surroundings to which heat is transferred by Newtonian cooling. 

The temperature dependence of the reactions is modeled with a temperature dependent reaction rate coefficient $k_0$. In fact all four reaction rate coefficients of (\ref{eq:reactSchema}) are temperature dependent. According to \cite{Scott1990} it
suffices, however, to consider only the temperature dependence of $k_0$. 

Following \cite{Scott1990}, we assume $P$ to be abundant and neglect its consumption. The resulting model (\ref{eq:ChemReactAB}) reads as 
\begin{eqnarray}\label{eq:dimensionlessReact}
\dot{\alpha} &=& \mu_0 e^{\delta \psi} - \alpha \beta^2 - \kappa_u \alpha, \nonumber\\
\dot{\beta} &=& \alpha \beta^2 + \kappa_u \alpha - \beta,\\
\dot{\psi} &=& \beta - \gamma \psi \nonumber 
\end{eqnarray}
with dimensionless concentrations $\alpha$ and $\beta$ of the chemical species $A$ and $B$, respectively, dimensionless temperature $\psi$, 
scaled initial concentration of $P$ 
$\mu_0=\sqrt{\frac{k_0^2 k_1}{k_2^3}} \, p_0,$  
rate constant $\kappa_u=k_3/k_2$, 
adiabatic temperature rise 
$\delta = (Q \sqrt{\frac{k_2}{k_1}}E)/(c_p c_0 R \tau^2_a)$,
activation energy $E$ 
of the first reaction in (\ref{eq:reactSchema}), 
ideal gas constant $R$, 
and the coefficient of Newtonian cooling 
$\gamma= (\chi S)/(k_2 V c_p c_0)$. 
The parameters $\kappa_u$ and $\delta$ are fixed to $\kappa_u=5.5\cdot 10^{-3}$ and $\delta=0.1$. 
The parameters $\mu_0$ and $\gamma$ are optimization variables.

We optimize \eqref{eq:dimensionlessReact} by maximizing the concentration of the final product~$C$.  Since the concentration of $C$ is proportional to the concentration of intermediate product $B$ in (\ref{eq:reactSchema}), we choose the objective function 
\begin{eqnarray}\label{eq:AverageConcentrForReact}
\phi=\frac{1}{T}\int_{0}^{T}\beta(t) \mathrm dt,
\end{eqnarray}
where $T$ is the period of the corresponding solution of system (\ref{eq:dimensionlessReact}). For equilibria this is equivalent to $\phi=\beta$. 

\subsubsection{Reference results obtained without normal vector constraints}
\begin{figure}
\centering
\includegraphics{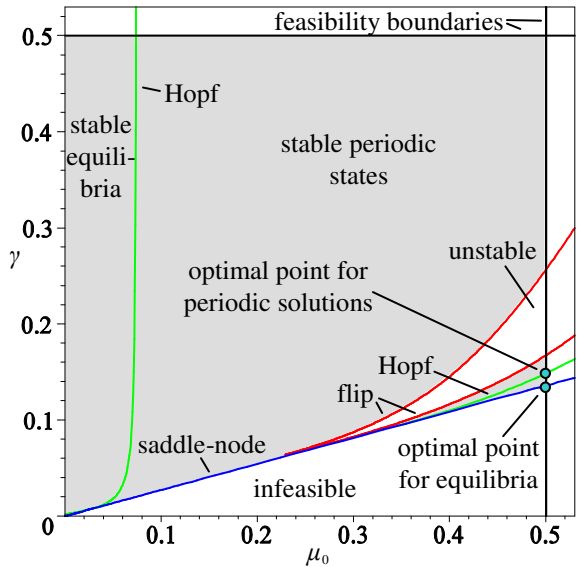}
\caption{Optimal equilibrium and periodic solutions for the reaction model~(\ref{eq:dimensionlessReact}) obtained without normal vector constraints. The robustness region \eqref{eq:UncRegReactionAB} is omitted in both cases for better visibility.}
\label{fig:ResABCReact}
\end{figure} 

We consider the feasibility constraints 
\begin{eqnarray}\label{eq:FeasibilityExample2}
  \mu_0 \leq 0.5, \quad \gamma \leq 0.5,   
\end{eqnarray}
the uncertainty region (\ref{eq:UncRegion}) 
\begin{eqnarray}\label{eq:UncRegReactionAB}
  (\mu_0,\gamma) \in ([\mu_0^{(0)}-\Delta \mu_0,\mu_0^{(0)}+\Delta \mu_0],[\gamma^{(0)}-\Delta \gamma,\gamma^{(0)}+\Delta \gamma]), 
\end{eqnarray}
where $\Delta \mu_0=\Delta \gamma=0.02$, and seek for the optimal equilibrium of system (\ref{eq:dimensionlessReact}).
Figure \ref{fig:ResABCReact} shows the optimal equilibrium that results from solving (\ref{eq:OptUnconstrainedSt}). The parameters and objective function evaluate to $(\mu_0^{(0)},\gamma^{(0)})=(0.5,0.1359)$ and $\phi=1.36$ at this point, respectively.

The optimal periodic solution that results from solving \eqref{eq:OptUnconstrainedPer}, which is also marked in Figure \ref{fig:ResABCReact}, corresponds to $(\mu_0^{(0)},\gamma^{(0)})=(0.5,0.1487)$ and $\phi=0.94$. This point results from optimizing over all periodic orbits without normal vector constraints, i.e., from solving \eqref{eq:OptUnconstrainedPer}.

The shaded areas in Figure \ref{fig:ResABCReact} correspond to stable and feasible modes of operation of the reaction system. Both optimal points  are located on the border of this area. Consequently, the optimal points that result from \eqref{eq:OptUnconstrainedPer} and \eqref{eq:OptUnconstrainedSt} are not robust, since there exist arbitrarily small parameter variations that result in a loss of stability.

We note for completeness that stable periodic solutions emanate from the Hopf bifurcations shown in Figure \ref{fig:ResABCReact}. Furthermore, stable periodic solutions lose stability due to flip bifurcations. This corroborates results by \cite{Scott1990}, who reported chaotic behavior as a result of cascades of period doubling.

\subsubsection{Results of the robust optimization with normal vector constraints}
\begin{figure}
\centering
\includegraphics{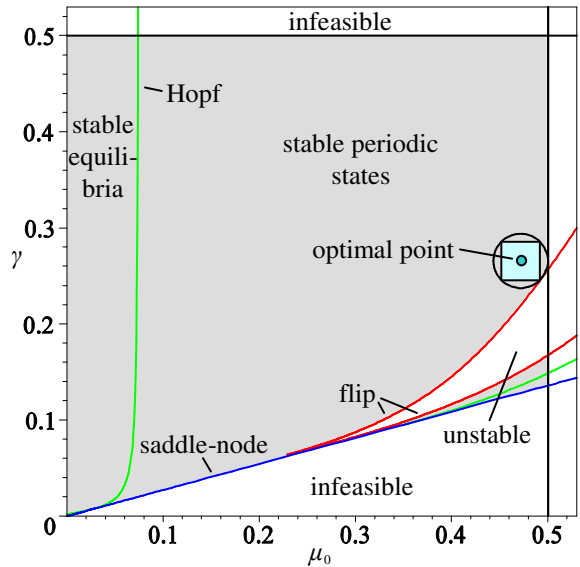}
\caption{Robust optimal point for the chemical reaction model~(\ref{eq:dimensionlessReact}) obtained with normal vector constraints.}
\label{fig:OptABCReact}
\end{figure}

\begin{figure}
\centering
\includegraphics{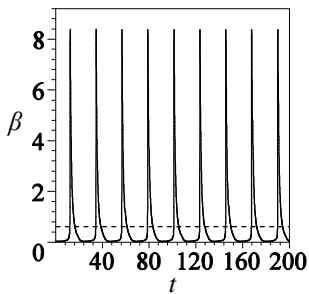}
\caption{Time series for the optimal point from Figure \ref{fig:OptABCReact}. The average concentration of $\beta$ (\ref{eq:AverageConcentrForReact}) is shown as the dashed line.}
\label{fig:OptTimeVsXABCReact}
\end{figure}

We solve the optimization problem \eqref{eq:OptProbelmWithNVConstr} to find the optimal stable and robust mode of operation. We initialize \eqref{eq:OptProbelmWithNVConstr} with a stable periodic solution.
All minimal distances to critical boundaries are set to $d_{\min}=\sqrt{n_{\alpha}}=\sqrt{2}$. 
Specifically, critical boundaries due to saddle-node bifurcation points of equilibria, flip bifurcation points of cycles, and the feasibility boundaries \eqref{eq:FeasibilityExample2} must be considered in this example. 
The resulting robust optimal point is illustrated in Figures \ref{fig:OptABCReact} and \ref{fig:OptTimeVsXABCReact}. It corresponds to $(\mu_0^{(0)},\gamma^{(0)})=(0.4717,0.2657)$ and the objective function value (\ref{eq:AverageConcentrForReact}) $\phi=0.61$. It is apparent from Figure \ref{fig:OptABCReact} that the optimal solution for \eqref{eq:OptProbelmWithNVConstr} is robust in the sense that there exists a stable and feasible solution for every combination of the uncertain parameters \eqref{eq:UncRegReactionAB}.
Note that it is not necessary here to switch between optimization problems (\ref{eq:OptProbelmWithNVConstr}) and (\ref{eq:OptProbelmWithNVConstrEquil}) in contrast to the previous example.

\section{Conclusion}\label{sec:Concl}
We extended the normal vector method for robust optimization of parametrically uncertain dynamical systems to the case of ODE systems with autonomous oscillations. It is the central idea of the proposed approach to use the Poincar\'e map to reduce the stability analysis of period orbits of continuous time systems to the stability analysis of fixed points of discrete time systems. By virtue of the Poincar\'e map, stability boundaries can be described with the bifurcation theory of fixed points. 
The proposed approach can naturally be combined with the normal vector method for equilibria whenever it is necessary to compare optimal and robust periodic orbits to optimal and robust equilibria.

We applied the proposed method to two chemical reaction processes that admit both robust equilibria and robust periodic orbits. A naive optimization that ignores stability properties yields an optimal mode of operation which, however, is unstable or not robust. We call an optimum not robust if there exists an arbitrarily small change of the optimal parameters that results in instability or infeasibility. In contrast to a naive optimization, the normal vector method provides an optimal robust mode of operation for both reaction systems.

\section*{Acknowledgment} Support by the Deutsche
    Forschungsgemeinschaft (DFG) under grant MO 1086/4 is gratefully acknowledged.

\appendix
\section{Augmented system for Neimark-Sacker bifurcations of cycles} \label{app:NSAugSys}
Necessary conditions for Neimark-Sacker bifurcation points of periodic solutions are given by the following set of $3n_x+3$ equations \citep{Lust1997}:
\begin{eqnarray}\label{eq:NSSystem}
M^{(\rm{NS})}(x_0,T,\alpha):=\left(\begin{array}{c}
\f(x_0,T,\alpha) -x_0\\
s(x_0,T,\alpha)\\
\f_{x_0}(x_0,T,\alpha)w^{(1)}-w^{(1)} \cos \theta +w^{(2)} \sin \theta\\
\f_{x_0}(x_0,T,\alpha)w^{(2)}-w^{(1)} \sin \theta -w^{(2)} \cos \theta\\
w^{(1) T}w^{(1)} + w^{(2) T}w^{(2)} - 1\\
w^{(1) T}w^{(2)}
\end{array}\right)=0.
\end{eqnarray}
The first two lines in (\ref{eq:NSSystem}) are the periodicity condition \eqref{eq:PerCondForPhi} and phase condition \eqref{eq:PhaseCondition}.
The third and forth line state that the Jacobian $\f_{x_0}$ has a pair of complex conjugate eigenvalues $e^{\pm i \theta}=\cos \theta \pm i \sin \theta$ corresponding to eigenvectors $w=w^{(1)} \pm i w^{(2)} \in \C^{n_x}$, respectively. The last two lines normalize the eigenvectors. The system of equations (\ref{eq:NSSystem}) is nonsingular with respect to $x_0$, $T$, $w^{(1)}$, $w^{(2)}$, $\theta$, and one component of $\alpha \in \R^{n_\alpha}$, say $\alpha_1$, at any nondegenerate Neimark-Sacker bifurcation point of cycles \citep{Lust1997}.

\section{Normal vector system for  Neimark-Sacker bifurcations of cycles} \label{app:NVSysNS}
In the case of Neimark-Sacker bifurcation points of periodic solutions the normal vector system reads as 
\begin{eqnarray}\nonumber
G^{(\rm{NS})}(p,\bar{x}^{(\rm{NS})},r):=\\
\left(\begin{array}{c}
M^{(\rm{NS})}(p)\\
\f_{x_0}^{T}(p) v^{(1)} - v^{(1)}  \cos \theta - v^{(2)} \sin \theta + \gamma_1 w^{(1)} - \gamma_2 w^{(2)}\\
\f_{x_0}^{T}(p) v^{(2)} + v^{(1)} \sin \theta - v^{(2)} \cos \theta + \gamma_1 w^{(2)} + \gamma_2 w^{(1)}\\
(w^{(1)T} v^{(1)} + w^{(2)T} v^{(2)} ) \sin \theta + (w^{(2)T} v^{(1)} - w^{(1)T} v^{(2)}) \cos \theta\\
v^{(1) T} w^{(1)} + v^{(2) T} w^{(2)} -1\\
\f_{x_0}^{T}(p) u - u + s_{x_0}^{T}(p) \varkappa +  v^{(1) T} \f_{x_0x_0}(p) w^{(1)} + v^{(2) T} \f_{x_0x_0}(p) w^{(2)}\\
\f_{T}^{T}(p) u + s_T(p) \varkappa +  v^{(1) T} \f_{x_0T}(p) w^{(1)} + v^{(2) T} \f_{x_0T}(p) w^{(2)}\\
r - \f_{\alpha}^{T}(p) u - s_{\alpha}^{T}(p) \varkappa -  v^{(1) T} \f_{x_0 \alpha}(p) w^{(1)} - v^{(2) T} \f_{x_0 \alpha}(p) w^{(2)} 
\end{array}\right)=0,\nonumber
\end{eqnarray}
where $p=(x_0,T,\alpha)$ as introduced in \eqref{eq:Abbreviations} and
\begin{eqnarray}\nonumber
\bar{x}^{(\rm{NS})}=(w^{(1)}, w^{(2)},\theta, v^{(1)}, v^{(2)}, u, \varkappa, \gamma_1, \gamma_2).
\end{eqnarray}
$M^{(\rm{NS})}(p)$ refer to system \eqref{eq:NSSystem}.
Vectors $w^{(1)}+iw^{(2)} \in \C^{n_x}$ and $v^{(1)}+iv^{(2)} \in \C^{n_x}$ are eigenvectors of the matrix $\f_{x_0}$ and its transpose $\f_{x_0}^T$ that correspond to the eigenvalues $e^{i \theta}$ and $e^{-i \theta}$, respectively, and $\gamma_2 \in \R$ is an auxiliary variable. 
All other symbols are defined as for the system $G^{(\rm{flip})}$ in Section~\ref{sec:NVDirect}.

\section{Normal vector system for Hopf and saddle-node bifurcations of equilibria}\label{app:NVSysSteady}
We state the normal vector systems for Hopf and saddle-node bifurcations of equilibria for ease of reference \citep{Moennigmann2002}. 
\begin{eqnarray}\nonumber
G^{(\rm{Hopf})}(q,\bar{x}^{(\rm{Hopf})},r):=\\
\left(\begin{array}{c}
f(q)\\
f_{x}(q)w^{(1)} + \omega w^{(2)}\\
f_{x}(q)w^{(2)} - \omega w^{(1)}\\
w^{(1)T}w^{(1)} + w^{(2)T}w^{(2)} - 1\\
w^{(1)T}w^{(2)}\\
f_{x}^{T}(q) v^{(1)} - \omega v^{(2)} + \gamma_1 w^{(1)} - \gamma_2 w^{(2)}\\
f_{x}^{T}(q) v^{(2)} + \omega v^{(1)} + \gamma_1 w^{(2)} + \gamma_2 w^{(1)}\\
v^{(1)T}w^{(1)} + v^{(2)T}w^{(2)} - 1\\
v^{(1)T}w^{(2)} - v^{(2)T}w^{(1)}\\
f_{x}^{T}(q) u +  v^{(1)T} f_{xx}(q) w^{(1)} +  v^{(2)T} f_{xx}(q) w^{(2)}\\
r - f_{\alpha}^{T}(q) u -  v^{(1)T} f_{x \alpha}(q) w^{(1)} -  v^{(2)T} f_{x \alpha}(q) w^{(2)}
\end{array}\right)=0,\nonumber
\end{eqnarray}
where $q=(x,\alpha)$
as in \eqref{eq:AbbreviationsII} and
$$\bar{x}^{(\rm{Hopf})}=(w^{(1)}, w^{(2)}, \omega, v^{(1)}, v^{(2)}, u, \gamma_1, \gamma_2).$$ 
Furthermore, 
$w^{(1)}+i w^{(2)} \in \C^{n_x}$ and $v^{(1)}+i v^{(2)} \in \C^{n_x}$ are eigenvectors of $f_x$ and $f_x^{T}$ corresponding to the eigenvalues $i \omega$ and $-i \omega$, respectively, $u \in \R^{n_x}$, $\gamma_1 \in \R$, and $\gamma_2 \in \R$ are auxiliary variables, and $r \in \R^{n_\alpha}$ denotes the normal vector.

The normal vector system for saddle-node bifurcations of equilibria reads
\begin{eqnarray}\nonumber
G^{(\rm{sn})}(q,\bar{x}^{(\rm{sn})},r):=
\left(\begin{array}{c}
f(q)\\
f_{x}^{T}(q) v\\
v^{T}v - 1\\
r - f_{\alpha}^{T}(q) v
\end{array}\right)=0,
\end{eqnarray}
where $q$ is as in \eqref{eq:AbbreviationsII} and $\bar{x}^{(\rm{sn})}$ equals $v$, the eigenvector of $f_{x}^{T}$ corresponding to eigenvalue zero. All other symbols are defined as for system $G^{(\rm{Hopf})}$ above. 

In contrast to the normal vector systems  $G^{(\rm{NS})}$ and  $G^{(\rm{flip})}$ for bifurcations of cycles, the derivatives $f_x$, $f_\alpha$, $f_{xx}$, and $f_{x \alpha}$ used for defining $G^{(\rm{Hopf})}$ and $G^{(\rm{sn})}$ can be obtained symbolically.

\bibliographystyle{elsarticle-harv}
\bibliography{dkabibl}

\begin{thebibliography}{44}
\expandafter\ifx\csname natexlab\endcsname\relax\def\natexlab#1{#1}\fi
\expandafter\ifx\csname url\endcsname\relax
  \def\url#1{\texttt{#1}}\fi
\expandafter\ifx\csname urlprefix\endcsname\relax\def\urlprefix{URL }\fi

\bibitem[{Abad et~al.(2009)Abad, Barrio, Blesa, and Rodr\'{i}guez}]{Abad2009}
Abad, A., Barrio, R., Blesa, F., Rodr\'{i}guez, M., 2009. {TIDES}: a {T}aylor
  {I}ntegrator for {D}ifferential {E}quation{S}. Preprint,
  http://gme.unizar.es/software/tides.

\bibitem[{Abashar and Elnashaie(2010)}]{Abashar2010}
Abashar, M. E.~E., Elnashaie, S. S. E.~H., 2010. Dynamic and chaotic behavior
  of periodically forced fermentors for bioethanol production. Chemical
  Engineering Science 65~(16), 4894--4905.

\bibitem[{Aguda et~al.(1989)Aguda, Larter, and Clarke}]{Aguda1989}
Aguda, B.~D., Larter, R., Clarke, B.~L., 1989. Dynamic elements of mixed-mode
  oscillations and chaos in a peroxidase-oxidase model network. Journal of
  Chemical Physics 90~(8), 4168--4175.

\bibitem[{Beyn et~al.(2002)Beyn, Champneys, Doedel, Govaerts, Kuznetsov, and
  Sandstede}]{Beyn2002}
Beyn, W.-J., Champneys, A., Doedel, E., Govaerts, W., Kuznetsov, Y.~A.,
  Sandstede, B., 2002. Numerical continuation, and computation of normal forms.
  In: Fiedler, B. (Ed.), Handbook of Dynamical Systems. Vol.~2. North--Holland,
  Amsterdam, pp. 149--219.

\bibitem[{Burke et~al.(2003)Burke, Lewis, and Overton}]{Burke2003}
Burke, J.~V., Lewis, A.~S., Overton, M.~L., 2003. Optimization and
  pseudospectra, with applications to robust stability. SIAM Journal on Matrix
  Analysis and Applications 25~(1), 80--104.

\bibitem[{Cash and Mazzia(2005)}]{Cash2005}
Cash, J.~R., Mazzia, F., 2005. A new mesh selection algorithm, based on
  conditioning, for two-point boundary value codes. Journal of Computational
  and Applied Mathematics 184~(2), 362--381.

\bibitem[{Chang and Sahinidis(2011)}]{Chang2011}
Chang, Y.~J., Sahinidis, N.~V., 2011. Steady-state process optimization with
  guaranteed robust stability under parametric uncertainty. AICHE Journal
  57~(12), 3395--3407.

\bibitem[{D'Avino et~al.(2006)D'Avino, Crescitelli, Maffettone, and
  Grosso}]{Avino2006}
D'Avino, G., Crescitelli, S., Maffettone, P.~L., Grosso, M., 2006. A critical
  appraisal of the {$\Pi$}-criterion through continuation/optimization.
  Chemical Engineering Science 61~(14), 4689--96.

\bibitem[{Diehl et~al.(2009)Diehl, Mombaur, and Noll}]{Diehl2009}
Diehl, M., Mombaur, K.~D., Noll, D., 2009. Stability optimization of hybrid
  periodic systems via a smooth criterion. IEEE Transactions on Automatic
  Control 54~(8), 1875--1880.

\bibitem[{Dobson(1993)}]{Dobson1993}
Dobson, I., 1993. Computing a closest bifurcation instability in
  multidimensional parameter space. Journal of Nonlinear Science 3~(3),
  307--327.

\bibitem[{Douglas and Rippin(1966)}]{Douglas1966}
Douglas, J.~M., Rippin, D. W.~T., 1966. Unsteady state process operation.
  Chemical Engineering Science 21~(4), 305--315.

\bibitem[{Engelborghs et~al.(1999)Engelborghs, Lust, and Roose}]{Lust1999}
Engelborghs, K., Lust, K., Roose, D., 1999. Direct computation of period
  doubling bifurcation points of large-scale systems of {ODE}s using a
  {N}ewton-{P}icard method. IMA Journal of Numerical Analysis 19~(4), 525--547.

\bibitem[{Gerhard et~al.(2008)Gerhard, Marquardt, and
  M\"{o}nnigmann}]{Gerhard2008}
Gerhard, J., Marquardt, W., M\"{o}nnigmann, M., 2008. Normal vectors on
  critical manifolds for robust design of transient processes in the presence
  of fast disturbances. SIAM Journal on Applied Dynamical Systems 7~(2),
  461--490.

\bibitem[{Gill et~al.(2001)Gill, Murray, Saunders, and Wright}]{Gill2001}
Gill, P.~E., Murray, W., Saunders, M.~A., Wright, M.~H., 2001. User's guide for
  {NPSOL} 5.0: A {F}ortran package for nonlinear programming. Systems
  Optimization Laboratory, Stanford University, Stanford, USA, {T}echnical
  report {SOL} 86-2.

\bibitem[{Halliwell(1978)}]{Halliwell1978}
Halliwell, B., 1978. Lignin synthesis - generation of hydrogen-peroxide and
  superoxide by horseradish-peroxidase and its stimulation by manganese {(II)}
  and phenols. Planta 140~(1), 81--88.

\bibitem[{Jianquiang and Ray(2000)}]{Jianquiang2000}
Jianquiang, S., Ray, A.~K., 2000. Performance improvement of activated sludge
  wastewater treatment by nonlinear natural oscillations. Chemical Engineering
  \& Technology 23~(12), 1115--1122.

\bibitem[{Kastsian(2012)}]{KastsianPhDThesis}
Kastsian, D., 2012. Robust optimization of discrete time systems and periodic
  operation with guaranteed stability. Ph.D. thesis, Ruhr-Universit\"{a}t
  Bochum.

\bibitem[{Kastsian and M\"{o}nnigmann(2010)}]{Kastsian2010}
Kastsian, D., M\"{o}nnigmann, M., 2010. Robust optimization of fixed points of
  nonlinear discrete time systems with uncertain parameters. SIAM Journal on
  Applied Dynamical Systems 9, 357--390.

\bibitem[{Kastsian and M\"{o}nnigmann(2012)}]{Kastsian2012}
Kastsian, D., M\"{o}nnigmann, M., 2012. Robust optimization of periodically
  operated reactors with stability constraints. In: Proceedings of the 2012
  IEEE International Conference on Control Applications. pp. 184--189.

\bibitem[{Khan et~al.(2007)Khan, Forouhar, Tao, and Tong}]{Khan2007}
Khan, J.~A., Forouhar, F., Tao, X., Tong, L., 2007. Nicotinamide adenine
  dinucleotide metabolism as an attractive target for drug discovery. Expert
  Opinion on Therapeutic Targets 11~(5), 695--705.

\bibitem[{Khinast and Luss(2000)}]{Khinast2000}
Khinast, J.~G., Luss, D., 2000. Efficient bifurcation analysis of
  periodically-forced distributed parameter systems. Computers \& Chemical
  Engineering 24~(1), 139--152.

\bibitem[{Kuznetsov(1998)}]{Kuznetsov1998}
Kuznetsov, Y.~A., 1998. Elements of applied bifurcation theory. Springer, New
  York.

\bibitem[{Larter(2003)}]{Larter2003}
Larter, R., 2003. Understanding complexity in biophysical chemistry. Journal of
  Physical Chemistry B 107~(2), 415--429.

\bibitem[{Lust(1997)}]{Lust1997}
Lust, K., 1997. Numerical bifurcation analysis of periodic solutions of partial
  differential equations. Ph.D. thesis, Katholieke Universiteit Leuven.

\bibitem[{M\"{a}der and F\"{u}ssl(1982)}]{Maeder1982}
M\"{a}der, M., F\"{u}ssl, R., 1982. Role of peroxidase in lignification of
  tobacco cells {II}. {R}egulation by phenolic compounds. Plant Physiology
  70~(4), 1132--1134.

\bibitem[{Mangold et~al.(2000)Mangold, Kienle, Gilles, and Mohl}]{Mangold2000}
Mangold, M., Kienle, A., Gilles, E.~D., Mohl, K.~D., 2000. Nonlinear
  computation in {DIVA} - methods and applications. Chemical Engineering
  Science 55~(2), 441--454.

\bibitem[{Mombaur(2009)}]{Mombaur2009}
Mombaur, K.~D., 2009. Using optimization to create self-stable human-like
  running. Robotica 27~(3), 321--330.

\bibitem[{Mombaur et~al.(2005{\natexlab{a}})Mombaur, Bock, Schloder, and
  Longman}]{Mombaur2005a}
Mombaur, K.~D., Bock, H.~G., Schloder, J.~P., Longman, R.~W.,
  2005{\natexlab{a}}. Open-loop stable solutions of periodic optimal control
  problems in robotics. ZAMM--Journal of Applied Mathematics and Mechanics
  85~(7), 499--515.

\bibitem[{Mombaur et~al.(2005{\natexlab{b}})Mombaur, Longman, Bock, and
  Schloder}]{Mombaur2005}
Mombaur, K.~D., Longman, R.~W., Bock, H.~G., Schloder, J.~P.,
  2005{\natexlab{b}}. Open-loop stable running. Robotica 23~(Part 1), 21--33.

\bibitem[{M\"{o}nnigmann and Marquardt(2002)}]{Moennigmann2002}
M\"{o}nnigmann, M., Marquardt, W., 2002. Normal vectors on manifolds of
  critical points for parametric robustness of equilibrium solutions of {ODE}
  systems. Journal of Nonlinear Science 12~(2), 85--112.

\bibitem[{M\"{o}nnigmann and Marquardt(2003)}]{Moennigmann2003}
M\"{o}nnigmann, M., Marquardt, W., 2003. Steady-state process optimization with
  guaranteed robust stability and feasibility. AICHE Journal 49~(12),
  3110--3126.

\bibitem[{M\"{o}nnigmann and Marquardt(2005)}]{Moennigmann2005}
M\"{o}nnigmann, M., Marquardt, W., 2005. Steady-state process optimization with
  guaranteed robust stability and flexibility: Application to {HDA} reaction
  section. Industrial \& Engineering Chemistry Research 44~(8), 2737--2753.

\bibitem[{M\"{o}nnigmann et~al.(2007)M\"{o}nnigmann, Marquardt, Bischof,
  Beelitz, Lang, and Willems}]{Moennigmann2007}
M\"{o}nnigmann, M., Marquardt, W., Bischof, C.~H., Beelitz, T., Lang, B.,
  Willems, P., 2007. A hybrid approach for efficient robust design of dynamic
  systems. SIAM Review 49~(2), 236--254.

\bibitem[{Mu\~{n}oz et~al.(2012)Mu\~{n}oz, Gerhard, and Marquardt}]{Munoz2012}
Mu\~{n}oz, D.~A., Gerhard, J., Marquardt, W., 2012. A normal vector approach
  for integrated process and control design with uncertain model parameters and
  disturbances. Computers \& Chemical Engineering 40, 202--212.

\bibitem[{Olsen(1983)}]{Olsen1993}
Olsen, L.~F., 1983. An enzyme reaction with a strange attractor. Physics
  Letters A 94~(9), 454--457.

\bibitem[{Parulekar(1998)}]{Parulekar1998}
Parulekar, S.~J., 1998. Analysis of forced periodic operations of continuous
  bioprocesses -- single input variations. Chemical Engineering Science
  53~(14), 2481--2502.

\bibitem[{Parulekar(2003)}]{Parulekar2003}
Parulekar, S.~J., 2003. Systematic performance analysis of continuous processes
  subject to multiple input cycling. Chemical Engineering Science 58~(23--24),
  5173--5194.

\bibitem[{Sauve(2008)}]{Sauve2008}
Sauve, A.~A., 2008. {NAD$^{+}$} and vitamin {B$_3$}: From metabolism to
  therapies. Journal of Pharmacology and Experimental Therapeutics 324~(3),
  883--893.

\bibitem[{Scott and Tomlin(1990)}]{Scott1990}
Scott, S.~K., Tomlin, A.~S., 1990. Period doubling and other complex
  bifurcations in non-isothermal chemical-systems. Philosophical Transactions
  of the Royal Society A 332~(1624), 51--68.

\bibitem[{Seydel(1988)}]{Seydel1988}
Seydel, R., 1988. From Equilibrium to Chaos. Practical Bifurcation and
  Stability Analysis. Elsevier, New York.

\bibitem[{Steinmetz et~al.(1993)Steinmetz, Geest, and Larter}]{Steinmetz1993}
Steinmetz, C.~G., Geest, T., Larter, R., 1993. Universality in the
  peroxidase-oxidase reaction: Period doublings, chaos, period three, and
  unstable limit cycles. Journal of Physical Chemistry 97~(21), 5649--5653.

\bibitem[{Sterman and Ydstie(1990)}]{Sterman1990}
Sterman, L.~E., Ydstie, B.~E., 1990. The steady-state process with periodic
  perturbations. Chemical Engineering Science 45~(3), 721--736.

\bibitem[{Stowers et~al.(2009)Stowers, Robertson, Ban, Tanner, and
  Boczko}]{Stowers2009}
Stowers, C.~C., Robertson, J.~B., Ban, H., Tanner, R.~D., Boczko, E.~M., 2009.
  Periodic fermentor yield and enhanced product enrichment from autonomous
  oscillations. Applied Biochemistry and Biotechnology 156~(1--3), 489--505.

\bibitem[{Vanbiervliet et~al.(2009)Vanbiervliet, Vandereycken, Michiels,
  Vandewalle, and Diehl}]{Vanbiervliet2009}
Vanbiervliet, J., Vandereycken, B., Michiels, W., Vandewalle, S., Diehl, M.,
  2009. The smoothed spectral abscissa for robust stability optimization. SIAM
  Journal on Optimization 20~(1), 156--171.

\end{thebibliography}

\end{document}